\documentclass[12pt]{article}
\usepackage{amsfonts,amsmath,amsxtra}
\usepackage{amssymb}
\def\hybrid{\topmargin 0pt      \oddsidemargin 0pt
     \headheight 0pt \headsep 0pt
     \textwidth 16.5cm
     \textheight 23cm
     \marginparwidth 0.0in
     \parskip 5pt plus 1pt   \jot = 1.5ex}
\catcode`\@=11
\def\marginnote#1{}
\newcount\hour
\newcount\minute
\newtoks\amorpm
\hour=\time\divide\hour by60
\minute=\time{\multiply\hour by60 \global\advance\minute by-\hour}
\edef\standardtime{{\ifnum\hour<12 \global\amorpm={am}%
     \else\global\amorpm={pm}\advance\hour by-12 \fi
     \ifnum\hour=0 \hour=12 \fi
   \number\hour:\ifnum\minute<10 0\fi\number\minute\the\amorpm}}
\edef\militarytime{\number\hour:\ifnum\minute<10 0\fi\number\minute}
\def\draftlabel#1{{\@bsphack\if@filesw {\let\thepage\relax
\xdef\@gtempa{\write\@auxout{\string
   \newlabel{#1}{{\@currentlabel}{\thepage}}}}}\@gtempa
\if@nobreak \ifvmode\nobreak\fi\fi\fi\@esphack}
     \gdef\@eqnlabel{#1}}
\def\@eqnlabel{}
\def\@vacuum{}
\def\draftmarginnote#1{\marginpar{\raggedright\scriptsize\tt#1}}
\def\draft{\oddsidemargin -0.1truein
     \def\@oddfoot{\sl preliminary draft \hfil
     \rm\thepage\hfil\sl\today\quad\militarytime}
     \let\@evenfoot\@oddfoot \overfullrule 3pt
     \let\label=\draftlabel
     \let\marginnote=\draftmarginnote
\def\@eqnnum{{\rm (\theequation)}
\rlap{\kern\marginparsep\tt\@eqnlabel}%
\global\let\@eqnlabel\@vacuum}  }
\newcommand{\RR}{{\mathbb{R}}}
\newcommand{\CC}{{\mathbb{C}}}
\newcommand{\NN}{{\mathbb{N}}}       

\newcommand{\ZZ}{{\mathbb{Z}}}

\newfont{\Bbbb}{msbm7 scaled 1\@ptsize00}
\newcommand{\zs}{\raise-1pt\hbox{$\mbox{\Bbbb Z}$}}

\font\sevenmsa=msam6 
\newfam\msafam
\textfont\msafam=\sevenmsa
\def\hexnumber@#1{\ifnum#1<10 \number#1\else
\ifnum#1=10 A\else\ifnum#1=11 B\else\ifnum#1=12 C\else
\ifnum#1=13 D\else\ifnum#1=14 E\else\ifnum#1=15 F\fi\fi\fi\fi\fi\fi\fi}
\def\msa@{\hexnumber@\msafam}
\def\llcorner{\delimiter"4\msa@78\msa@78 }
\def\lrcorner{\delimiter"5\msa@79\msa@79 }
\mathchardef\blacktriangleright="3\msa@49
\mathchardef\blacktriangleleft="3\msa@4A
\font\tenmsb=msbm10 scaled 1\@ptsize00
\newfam\msbfam
\textfont\msbfam=\tenmsb
\scriptfont\msbfam=\tenmsb
\def\msb@{\hexnumber@\msbfam}
\mathchardef\varkappa="0\msb@7B
\newdimen\linethick  \linethick=0.4pt
\newdimen\hboxitspace    \hboxitspace=5pt
\newdimen\vboxitspace    \vboxitspace=5pt
\def\fr#1{%
\be\new
\vcenter{
\hrule height\linethick
        \hbox{\vrule width\linethick
              \kern\hboxitspace
              \vbox{\kern\vboxitspace
                    \hbox{$\begin{array}{c}\displaystyle#1
       \end{array}$}%
                    \kern\vboxitspace}%
              \kern\hboxitspace
              \vrule width\linethick}%
        \hrule height\linethick}%
\ee}
\newdimen\Squaresize \Squaresize=14pt
\newdimen\Thickness \Thickness=0.5pt
\def\Square#1{\hbox{\vrule width \Thickness
\vbox to \Squaresize{\hrule height \Thickness\vss
   \hbox to \Squaresize{\hss#1\hss}
\vss\hrule height\Thickness}
\unskip\vrule width \Thickness}
\kern-\Thickness}
\def\Vsquare#1{\vbox{\Square{$#1$}}\kern-\Thickness}

\def\numberbysection{\@addtoreset{equation}{section}
     \def\theequation{\thesection.\arabic{equation}}}
\numberbysection

\renewcommand{\theequation}{\thesection.\arabic{equation}}
\def\titlepage{\@restonecolfalse\if@twocolumn\@restonecoltrue\onecolumn
  \else \newpage \fi \thispagestyle{empty}\c@page\z@
     \def\thefootnote{\fnsymbol{footnote}} }
\def\endtitlepage{\if@restonecol\twocolumn \else  \fi
     \def\thefootnote{\arabic{footnote}}
     \setcounter{footnote}{0}}  
\relax
\hybrid
\makeatletter
\newdimen\normalarrayskip            
\newdimen\minarrayskip               
\normalarrayskip\baselineskip
\minarrayskip\jot
\newif\ifold             \oldtrue            \def\new{\oldfalse}
\def\arraymode{\ifold\relax\else\displaystyle\fi}
\def\eqnumphantom{\phantom{(\theequation)}} 
\def\@arrayskip{\ifold\baselineskip\z@\lineskip\z@
  \else
  \baselineskip\minarrayskip\lineskip1\baselineskip\fi}
\def\@arrayclassz{\ifcase \@lastchclass \@acolampacol \or
\@ampacol \or \or \or \@addamp \or
\@acolampacol \or \@firstampfalse \@acol \fi
\edef\@preamble{\@preamble
\ifcase \@chnum
  \hfil$\relax\arraymode\@sharp$\hfil
  \or $\relax\arraymode\@sharp$\hfil
  \or \hfil$\relax\arraymode\@sharp$\fi}}
\def\@array[#1]#2{\setbox\@arstrutbox=\hbox{\vrule
  height\arraystretch \ht\strutbox
  depth\arraystretch \dp\strutbox
width\z@}\@mkpream{#2}\edef\@preamble{\halign \noexpand\@halignto
\bgroup \tabskip\z@ \@arstrut \@preamble \tabskip\z@ \cr}%
\let\@startpbox\@@startpbox \let\@endpbox\@@endpbox
\if #1t\vtop \else \if#1b\vbox \else \vcenter \fi\fi
\bgroup \let\par\relax
\let\@sharp##\let\protect\relax
\@arrayskip\@preamble}
\def\eqnarray{\stepcounter{equation}%
           \let\@currentlabel=\theequation
           \global\@eqnswtrue
           \global\@eqcnt\z@
           \tabskip\@centering              
           \let\\=\@eqncr
           $$%
         \halign to \displaywidth  \bgroup
          \eqnumphantom \@eqnsel
   \hskip\@centering                               
 $\displaystyle  \tabskip\z@ {##}$%
 &\global\@eqcnt\@ne \hskip 2\arraycolsep
      $ \displaystyle  \arraymode{##}$\hfil
 &\global\@eqcnt\tw@ \hskip 2\arraycolsep
      $\displaystyle\tabskip\z@{##}$\hfil
      \tabskip\@centering
 &{##}\tabskip\z@\cr}
\makeatother

\newtheorem{prop}{Proposition}[section]           
\newtheorem{cor}{Corollary}[section]
\newtheorem{lem}{Lemma}[section]
\newtheorem{ex}{Example}[section]
\newtheorem{rem}{Remark}[section]
\newcommand{\beq}[1]{\begin{equation}\label{#1}}
\newcommand\eeq{\end{equation}}
\newcommand\bqa{\begin{eqnarray}}
\newcommand\eqa{\end{eqnarray}}
\def\be{\begin{eqnarray}\new\begin{array}{cc}}
\def\ee{\end{array}\end{eqnarray}}

\def\ws{\hfill{$\square$}}                              
\def\beq{\begin{equation}}
\def\eeq{\end{equation}}
\def\bse{\begin{subequations}}                
\def\ese{\end{subequations}}
\def\bp{\begin{pmatrix}}
\def\ep{\end{pmatrix}}
\def\bs{\begin{subequations}\label}
\def\es{\end{subequations}}
\def\beq{\begin{eqnarray}{cc}}
\def\eeq{\end{eqnarray}}
\newcommand{\ba}{\begin{align}}
\newcommand{\ea}{\end{align}}
\def\bp{\begin{pmatrix}}
\def\ep{\end{pmatrix}}
\def\bv{\begin{vmatrix}}
\def\ev{\end{vmatrix}}
\def\bel{\be\label}

\def\stack#1#2{\raise0.7pt\hbox{$\mathrel{\mathop{#2}\limits^{#1}}$}}
\def\tr{\triangleright}
\def\tl{\triangleleft}
\def\sem{\mathsurround=0pt \raise1pt
\hbox{$\scriptscriptstyle>\!\!$}\:\!\!\tl}
\def\mes{\mathsurround=0pt \tr\!\:\!\raise0.8pt
\hbox{$\scriptscriptstyle\!\!<$}\,}
\def\]{\mathsurround=0pt ]\raise-2pt\hbox{$_\ast$}}

\def\<{\langle}
\def\>{\rangle}

\def\we{\raise-1pt\hbox{$\,\stackrel{\wedge}{,}\,$}}

\begin{document}
\setcounter{footnote}0
\title{
\hfill{\normalsize ITEP-TH-22/15}\\
\hfill{\normalsize IITP/TH-14/15}
\\ [10mm] \bf  Theta vocabulary II.\\
Multidimensional case}
\author{S. Kharchev\thanks{Institute for Theoretical and Experimental
Physics, Moscow, Russia and Institute for Information Transmission Problems, Moscow 127994, Russia}
\and A. Zabrodin\thanks{National Research University
Higher School of Economics,
International Laboratory of Representation Theory
and Mathematical Physics,
20 Myasnitskaya Ulitsa, Moscow 101000, Russia and
Institute for Theoretical and Experimental
Physics, Moscow, Russia}}
\date{\phantom{.}}
\maketitle

\begin{abstract}
It is shown that the Jacobi and Riemann identities of degree four
for the multidimensional theta functions
as well as the Weierstrass identities emerge as algebraic
consequences of the fundamental multidimensional binary identities
connecting the theta functions with Riemann matrices $\tau$ and $2\tau$.
\end{abstract}

\clearpage \newpage

\section{Introduction}

In the previous paper \cite{KZ}, we have shown that all known identities for Jacobi theta functions of one variable with half-integer
characteristics (of the Jacobi, Weierstrass and Riemann type) are algebraic
corollaries of the six fundamental 3-term binary identities that connect theta functions
with modular parameters $\tau$ and $2\tau$. In this paper we give a
multidimensional extension of these results.

Let $g\in\NN$, $\<.\,,.\>:\,\CC^g\times\CC^g\to\CC$ be the standard scalar product. Let $u:=(u_1,\ldots u_g)\in\CC^g$, and let $\tau$ be a symmetric $g\times g$ matrix
with positively definite imaginary part.
For any $a,b\in\RR^g$ the $g$-dimensional theta functions with characteristics
are defined as follows:
\bel{in0}
\theta_{[a;b]}(u|\tau)=
\sum_{k\in\ZZ^g}\exp\{\pi i\<\tau(k+a),k+a\>+2\pi i\<k+a,u+b\>\}.
\ee
The claim is that the binary relations (due to Schr\"oter \cite{Sch})
\bel{in1}
\theta_{[a_1;b_1]}(u_1|\tau)\theta_{[a_2; b_2]}(u_2|\tau)\\=
\sum_{p\in\ZZ^g/2\ZZ^g}
\theta_{[\frac{a_1+a_2+p}{2}; b_1+b_2]}(u_1+u_2|2\tau)
\theta_{[\frac{a_1-a_2+p}{2}; b_1-b_2]}(u_1-u_2|2\tau)
\ee
are fundamental ones in the sense that all
4th order identities (Jacobi and Riemann) and Weierstrass identities of order $2^g+2$
can be derived from the Schr\"oter relations. We also show that the
identities of the Jacobi and Riemann types are equivalent.

To formulate the main statement, it is necessary to introduce some notations.
Let $x_k\in\CC^g,\,k=1,2,3,4$ be $g$-dimensional vectors. Define the Whittaker-Watson dual vectors $x'_k$ as follows \cite{WW}:
\bel{in3}
\ \ x'_1={\textstyle\frac{1}{2}}(-x_1+x_2+x_3+x_4),\\
x'_2={\textstyle\frac{1}{2}}(x_1-x_2+x_3+x_4),\\
x'_3={\textstyle\frac{1}{2}}(x_1+x_2-x_3+x_4),\\
x'_4={\textstyle\frac{1}{2}}(x_1+x_2+x_3-x_4).
\ee
\begin{prop}\label{prop1}
\begin{itemize}
\item[(i)] As a consequence of (\ref{in1}), the following
multidimensional Jacobi identities hold:
\bel{in5}
\hspace{-1cm}
\sum_{q\in\ZZ^g/2\ZZ^g}\prod_{k=1}^4
e^{-\pi i\<a_k,q\>}\theta_{[a_k;b_k+\frac{q}{2}]}(u_k|\tau)=
\sum_{q\in\ZZ^g/2\ZZ^g}\prod_{k=1}^4
e^{-\pi i\<a_k',q\>}\theta_{[a_k';b_k'+\frac{q}{2}]}(u_k'|\tau),
\ee
where the dual variables are defined by (\ref{in3}).

\smallskip\noindent
\item[(ii)] Similarly, there are Jacobi identities:
\bel{in6}
\sum_{p\in\ZZ^g/2\ZZ^g}
\prod_{k=1}^4\theta_{[a_k+\frac{p}{2};b_k]}(u_k|\tau)=
\sum_{p\in\ZZ^g/2\ZZ^g}
\prod_{k=1}^4\theta_{[a_k'+\frac{p}{2};b_k']}(u_k'|\tau).
\ee
\item[(iii)] As a consequence of (\ref{in5}) or (\ref{in6}), there
are multidimensional Riemann identities:
\bel{in7}
\hspace{-0.3cm}
\prod_{k=1}^4\theta_{[a_k';b_k']}(u_k'|\tau)=2^{-g}
\sum_{p,q\in\ZZ^g/2\ZZ^g}e^{-\pi i\<p,q\>}\prod_{k=1}^4
e^{-\pi i\<a_k,q\>}\theta_{[a_k+\frac{p}{2};b_k+\frac{q}{2}]}(u_k|\tau).
\ee
\item[(iv)]
The identities (\ref{in5}), (\ref{in6}), and (\ref{in7}) are equivalent.
\item[(v)]
Let $\theta(u|\tau)$ be an arbitrary $g$-dimensional odd theta function. Then, by virtue of (\ref{in1}), the following Weierstrass identity holds:
\bel{pf1}
{\rm Pf}||\theta(w_i+w_j)\theta(w_i-w_j)||_{i,j=1}^{2^g+2}=0,
\ee
where ${\rm Pf}||a_{i,j}||$ is the Pfaffian of the skew-symmetric matrix $||a_{i,j}||$.
\end{itemize}
\end{prop}

Let us make some historical remarks. Among the huge diversity of non-trivial identities satisfied by the theta functions, the binary identities are the oldest ones. The history can be traced to Jacobi's paper \cite{J1} published in 1828. Namely, Jacobi wrote the identity $H(x,q)\Theta(X,q)-H(X,q)\Theta(x,q)=H(\frac{x-X}{2},\sqrt q)
H(\frac{\pi}{2}-\frac{x+X}{2},\sqrt q)$
\footnote{
 The Jacobi functions $H(x,q)$ and $\Theta(x,q)$ (where $q:=e^{i\pi\tau}$) 
 correspond to the standard functions $\theta_1(x|\tau)$ and $\theta_4(x|\tau)$, respectively.
} and noted that 
{\it l'\'equation remarquable et ais\'ement \`a d\'emontrer au moyen des premiers \'el\'ements de la trigonom\'etrie}
\cite[p. 305]{J1}. The general one-dimensional binary identities are due to Schr\"oter
\cite{Sch} (1854) who wrote the formulas for the product of theta functions with modular parameters $n_1\tau$ and $n_2\tau,\,(n_1,n_2\in\NN)$ (see below). 
One-dimensional four- and five-term identities of degree 4 were derived by Jacobi \cite{J2} in 1835-1836. The three-term identity of degree 4 (the addition formula) have been obtained by Weierstrass \cite{We2} in 1862. It  turns out that all one-dimensional identities of degree 4 are equivalent (see \cite{K} and \cite{KZ}).

Multidimensional binary identities have been written by K\"onigsberger in 1864 (without proof) \cite{Ko}.
A generalization of 4th order five-term Jacobi identities was first given by Prym  \cite{Pr}. Prym named the multidimensional identities {\it the Riemann theta formulas}.
In the preface of his book, Prym wrote that he learned of the formulas during his meeting with Riemann in Pisa in early 1865, and stressed that he wrote down a proof following
Riemann's suggestions
\footnote{
Im Fruhjahre 1865 war mir das Gluck zu Theil geworden, bei meinem hoch verehrten Lehrer Riemann in Pisa, wo derselbe sich seiner Gesundheit wegen aufhielt, einige Wochen zubringen zu konnen. Ich war damals mit gewissen Untersuchungen aus der Theorie der hyperelliptischen Functionen beschaftig, deren Anfange ich schon in meiner in Jahre 1864 den Denkschriften der Weiner Akademie erschienenen Arbeit "Neue Theorie der ultraelliptishen Functionen" anhangsweise veroffentlicht hatte, und es bildete wahrend meins Aufenthaltes in Pisa unter anderem auch das Additiontheorem der hyperelliptishen Functionen einen Gegenstand meiner Studien. Bei dieser Gelegenheit wurde mir von Riemann eine Formel (Formel (12) der esten Arbeit) mitgetheilt, die fur die Theorie der Thetafunctionen als eine fundamentale anzusehen ist, und ich verfasste auf seine Anregung hin einen Beweis fur diese Formel, dessen Gang auch die Zustimmung meines Lehrers fand. Zu einer Verwerthung der erwahnten Formel gelangte ich aber damals nicht, einmal, weil eine Verschlimmerung in dem Befinden Riemann's weitere Besprechungen unmoglich machte, dann aber auch, weil eingehendere auf die Charakteristiken der Thetareihen bezugliche  Untersuchungen, deren vorherige Durchfuhrung mir nothwendig erschien, mich ganz in Anspruch nahmen. \cite[p. V]{Pr}
}.
In \cite{We1}, Weierstrass presented the multidimensional generalizations
of both four-term and five-term identities of the 4th order in his dissertation which has not been published at that time ('due to typographical difficulties'). He also mentioned that the Prym's results are the particular case of his formulas
\footnote{
Specielle Falle der im Vorstehenden entwickelten Gleichungen (I.), (II.), (III.), (IV.) sind bereits von Anderen , namentlich von Herrn Prym, behandelt worden. Ich habe meine Abhandlung bereits vor einer Reihe von Jahren der Akademie vorgelegt, dieselbe konnte jedoch damals wegen typographischer Schwierigkeiten nicht veroffentlicht werden. Ich mochte aber noch bemerken, dass die Gleichungen (IL), (III.), (IV.) aus der Gleichung (I.) auch dadurch erhalten werden konnen, dass man die darin vorkommenden Argumente um gewisse Constanten (halbe Perioden der zu den vorkommenden $\Theta$-Functionen gehorenden Abelschen Integrale erster Art) vermehrt, wie dies in dem Falle n = l von Jacobi geschehen ist, wahrend bei meiner Ableitung die Kenntniss jener Constanten nicht erforderlich ist. \cite[p. 137]{We1}.
}. Quite general multidimensional identities can be found in \cite{M1}, \cite{M2}.

The generalization of the
addition formula to multidimensional case in terms of Pfaffian is due to Weierstrass \cite{We2}.

The main goal of the paper is to represent the detailed proof of Proposition \ref{prop1} thus establishing the fact that the binary relations are fundamental ones to compare with all other higher order identities.

In section 2, we give the definition of the multidimensional theta functions and discuss their properties which are essential for the proof.

In section 3, we represent the detailed derivation of the general multidimensional Schr\"oter identities which describe the product
$\theta_{[a_1;b_1]}(u_1|n_1\tau)\theta_{[a_2;b_2]}(u_2|n_2\tau)$ as a linear combination of appropriate products of theta functions with Riemann matrices
$(n_1+n_2)\tau$ and $n_1n_2(n_1+n_2)\tau$, where $n_1,n_2\in\NN$.

In section 4 and 5, we prove that the Jacobi identities of both types (\ref{in5}), (\ref{in6}) as well as the Riemann identities (\ref{in7}), are the simple corollaries of the standard Schr\"oter relations ($n_1=n_2=1$). For completeness, we represent all identities both in terms of Whittaker-Watson variables $u_k'$ defined by (\ref{in3}) and Jacobi ones, $\tilde u_k$ \cite[p.503]{J2} (see definition (\ref{J7}) below) which are also widely used (see \cite[p. 212]{M1}, \cite[p.102]{M2}).

The equivalence of identities (\ref{in5}), (\ref{in6}), and (\ref{in7}) is proved in section 6. It is shown that there are equivalent relations (see (\ref{R7}) below) which can be considered as a certain generalization of Weierstrass addition formulas (see complete list for $g=1$ in \cite{KZ}). These {\it naive} Weierstrass identities of degree 4 relate the products of theta function written in terms three sets of variables $u_k,u_k',\tilde u_k, (k=1,2,3,4)$.

Finally, in section 7, we present the simple proof of original Weierstrass identities (\ref{pf1}) of degree $2^g+2$ \cite{We2} with the help of multidimensional Shr\"oter's relations (\ref{in1}).

\paragraph{Acknowledgments.}
This work was supported in part by grant NSh-1500.2014.2
for support of scientific schools. The work of \,S.K. was supported in part by RFBR grant 15-01-99504. The work of A.Z. was supported in part by RFBR grant
14-02-00627 and joint grants 15-52-50041-YaF, 14-01-90405-Ukr.
The article was prepared within the framework of a subsidy granted to the HSE by the Government of the Russian Federation for the implementation of the Global Competitiveness Program.

\section{Multidimensional theta functions}
Let $g\in\NN$, $\<.\,,.\>:\,\CC^g\times\CC^g\to\CC$ be the standard scalar product. Let $u:=(u_1,\ldots u_g)\in\CC^g$, and symmetric $g\times g$ matrix
$\tau$ has positive definite imaginary part. Consider the infinite series \cite{M1}:
\bel{th}
\theta_{[a;b]}(u|\tau)=
\sum_{k\in\ZZ^g}\exp\{\pi i\<\tau(k+a),k+a\>+2\pi i\<k+a,u+b\>\},
\ee
where $i=\sqrt{-1}$ and $a,b\in\RR^g$. The series is absolutely convergent for
any $u\in \CC$ and defines the entire function $\theta_{[a;b]}(u|\tau)$.
It is called $g$-dimensional theta function with characteristics $[a;b]$.

It directly follows from definition that the following relations hold:
\bs{d3}
\ba
\theta_{[a+m;b+n]}(u|\tau)&=e^{2\pi i\<a,n\>}\theta_{[a;b]}(u|\tau),\ \ \
m,n\in\ZZ^g,\label{d3a}\\
\theta_{[-a;-b]}(-u|\tau)&=\theta_{[a;b]}(u|\tau).\label{d3b}
\end{align}
\es
By virtue of (\ref{d3}), it is sufficient
to consider the theta functions with characteristics $[a;b]$ such that
$0\leq a_j,b_j<1,\,j=1,\ldots,g$.

The characteristics $[a,b]$ for which all components $a_j,b_j$ are 0 or $\frac{1}{2}$ are called {\it half-periods}. A half-period $[a,b]$ is said to be {\it even} if
$4\<a,b\>= 0\,({\rm mod}\,2)$ and {\it odd} otherwise.
From (\ref{d3}) it follows that
$\theta_{[a;b]}(-u|\tau)=e^{4\pi i\<a,b\>}\theta_{[a;b]}(u|\tau),\,
(a,b\in\frac{1}{2}\ZZ^g/2\ZZ^g)$,
i.e. the function $\theta_{[a;b]}(u|\tau)$ is even or odd according to whether $[a;b]$ is even or odd half-periods.
It is well-known (see \cite{M1} for example) that there are $2^{g-1}(2^g+1)$
even half-periods and $2^{g-1}(2^g-1)$ odd ones.

Note that the functions (\ref{th}) with different characteristics are connected by the relations
\bel{d4}
\theta_{[a;b]}(u+\tau a'+b'|\tau)=
e^{-\pi i\<\tau a',a'\>-2\pi i\<a',u+b+b'\>}\theta_{[a+a';b+b']}(u|\tau)
\ee
for all $a',b'\in\RR^g$. In particular, by virtue of (\ref{d3a}),
\be
\theta_{[a;b]}(u+n|\tau)=e^{2\pi i\<a,n\>}\theta_{[a;b]}(u|\tau),\\
\theta_{[a;b]}(u+\tau n|\tau)=e^{-\pi i\<n,2u+2b+\tau n\>}\theta_{[a;b]}(u|\tau).
\ee
Hence, the functions $\theta_{[a;b]}(u|\tau)$ are quasiperiodic with (quasi)periods
$n$ and $\tau n,\;n\in\ZZ^g$.

There is a plethora of relations satisfied by the theta functions.
The simplest ones are well-known linear identities which relate the
theta functions with the Riemann matrices $\tau$ and $n^2\tau, n\in\NN$.
\begin{prop}
Let $n\in\NN$. The following identities hold:
\bs{sr1}
\ba
\theta_{[a;b]}(u|\tau)&=\sum_{p\in\ZZ^g/n\ZZ^g}\theta_{[\frac{a+p}{n};nb]}(nu|n^2\tau),
\label{sr1a}\\
\theta_{[a;b]}(nu|n^2\tau)&=n^{-g}\sum_{q\in\ZZ^g/n\ZZ^g}
e^{-2\pi i\<a,q\>}\theta_{[na;\frac{b+q}{n}]}(u|\tau).\label{sr1b}
\end{align}
\es
\end{prop}
{\bf Proof.} Starting with definition (\ref{th}), one can
represent the summation over $k\in\ZZ^g$ as follows:
$k=nl+p,\,l\in\ZZ^g,\,p\in\ZZ^g/n\ZZ^g$. Then
$\sum_{k\in\ZZ}(\ldots)=\sum_{p\in\ZZ^g/n\ZZ^g}\sum_{l\in\ZZ}(\ldots)$ and (\ref{sr1a}) is proved. To prove (\ref{sr1b}), perform the shift $b\to b+\frac{q}{n},\,q\in\ZZ^g/n\ZZ^g$ in (\ref{sr1a}). By virtue of (\ref{d3a}), one obtains the relation
\bel{sr2}
e^{-\frac{2\pi i}{n}\<a,q\>}\theta_{[a;b+\frac{q}{n}]}(u|\tau)=
\sum_{p\in\ZZ^g/n\ZZ^g}e^{\frac{2\pi i}{n}\<p,q\>}\theta_{[\frac{a+p}{n};nb]}(nu|n^2\tau).
\ee
Performing the summation over $q\in\ZZ^g/n\ZZ^g$ and using the identity
\bel{sr3}
\sum_{q\in\ZZ^g/n\ZZ^g}e^{\pm\frac{2\pi i}{n}\<p,q\>}=n^g\delta_{p,0},
\ee
which holds for any $p\in\ZZ^g/n\ZZ^g$, one arrives to (\ref{sr1b}). \ws

Below, we derive more complicated identities of higher order.

\section{Schr\"oter binary identities}
In this section, we represent the simple derivation of multidimensional binary relations which are the building blocks for higher order identities. The material of this section is not new, but we present all the details for completeness.
\begin{prop}{\rm \cite{Sch}}.
Let $n_1,n_2\in\NN$. The following binary identities hold:
\bel{Sch1}
\theta_{[a_1;b_1]}(u_1|n_1\tau)\theta_{[a_2;b_2]}(u_2|n_2\tau)\\
=\sum_{p\in\ZZ^g/(n_1+n_2)\ZZ^g}
\theta_{[\frac{n_1a_1+n_2a_2+n_1p}{n_1+n_2};b_1+b_2]}(u_1+u_2|(n_1+n_2)\tau)\\ \times\,
\theta_{[\frac{a_1-a_2+p}{n_1+n_2};n_2b_1-n_1b_2]}
(n_2u_1-n_1u_2|n_1n_2(n_1+n_2)\tau).
\ee
\end{prop}
{\bf Proof.} The proof is essentially the one given in \cite{Sch}.
Consider the simplest case of zero characteristics $a_k=b_k=0,\,k=1,2$.
One has the product  $\theta_{[0;0]}(u_1|n_1\tau)\theta_{[0;0]}(u_2|n_2\tau)=\sum_{k_1,k_2\in\ZZ^g}
e^{\pi i F_2+2\pi iF_1}$, where $F_2:=n_1\<\tau k_1,k_1\>+n_2\<\tau k_2,k_2\>$, $F_1:=\<k_1,u_1\>+\<k_2,u_2\>$. First we transform the quadratic part $F_2$.
The whole procedure consists of three simple steps.

Changing
$k_1:=l_1+k_2,l_2\in\ZZ^g$, we have:
$F_2=n_1\<\tau(l_1+k_2),l_1+k_2\>+n_2\<\tau k_2,k_2\>$.
This can be identically written as
$F_2=\frac{n_1n_2}{n_1+n_2}\<\tau l_1,l_1\>+(n_1+n_2)\<\tau(k_2+\frac{n_1}{n_1+n_2}l_1),k_2+\frac{n_1}{n_1+n_2}l_1\>$.

Next, one can divide the sum over $l_1$ into two sums setting $l_1:=(n_1+n_2)m_1+p$, where $m_1\in\ZZ^g$, $p\in\ZZ^g/(n_1+n_2)\ZZ^g$.
Thus $F_2=n_1n_2(n_1+n_2)
\<\tau(m_1+\frac{p}{n_1+n_2}),m_1+\frac{p}{n_1+n_2}\>+$ $(n_1+n_2)\<\tau(k_2+n_1m_1+\frac{n_1p}{n_1+n_2}),k_2+n_1m_1+\frac{n_1p}{n_1+n_2}\>$.

Changing the summation index $k_2:=m_2-n_1m_1$, $m_2\in\ZZ^g$, the term $F_2$ acquires the final form:
$F_2=n_1n_2(n_1+n_2)\<\tau(m_1+\frac{p}{n_1+n_2}),m_1+\frac{p}{n_1+n_2}\>+$
$(n_1+n_2)\<\tau(m_2+\frac{n_1p}{n_1+n_2}),m_2+\frac{n_1p}{n_1+n_2}\>$, where $m_1,m_2\in\ZZ^g, p\in\ZZ^g/(n_1+n_2)\ZZ^g$.

Evidently, after the same changes of summation indices, the linear part acquires the form:
$F_1=\<m_1+\frac{p}{n_1+n_2},n_2u_1-n_1u_2\>+\<m_2+\frac{n_1p}{n_1+n_2},u_1+u_2\>$.
Thus one proves the identity
\be
\theta_{[0;0]}(u_1|n_1\tau)\theta_{[0;0]}(u_2|n_2\tau)\\ \hspace{-1cm}=
\sum_{p\in\ZZ^g/(n_1+n_2)\ZZ^g}
\theta_{[\frac{n_1p}{n_1+n_2};0]}(u_1+u_2|(n_1+n_2)\tau)
\theta_{[\frac{p}{n_1+n_2};0]}(n_2u_1-n_1u_2|n_1n_2(n_1+n_2)\tau).\hspace{-0.5cm}
\ee
Using the relation (\ref{d4}), one arrives to the general formula (\ref{Sch1}). \ws

\smallskip
Consider the particular case $n_1=n_2:=n$. Changing $n\tau\to\tau$, one arrives to the
relation:
\bel{Sch2}
\theta_{[a_1;b_1]}(u_1|\tau)\theta_{[a_2;b_2]}(u_2|\tau)\\
=\sum_{p\in\ZZ^g/2n\ZZ^g}\theta_{[\frac{a_1+a_2+p}{2};b_1+b_2]}(u_1+u_2|2\tau)
\theta_{[\frac{a_1-a_2+p}{2n};n(b_1-b_2)]}(n(u_1-u_2)|2n^2\tau).
\ee
Let us show in detail that the relation (\ref{Sch2}) can be reduced to the standard form with $n=1$.
Indeed, using the linear identity (\ref{sr1b}), one can write
$\theta_{[\frac{a_1-a_2+p}{2n};n(b_1-b_2)]}(n(u_1-u_2)|2n^2\tau)=n^{-g}
\sum_{q\in\ZZ^g/n\ZZ^g}e^{-\frac{\pi i}{n}\<a_1-a_2+p,q\>}
\theta_{[\frac{a_1-a_2+p}{2};b_1-b_2+\frac{q}{n}]}(u_1-u_2|\tau)$. Substituting this identity to (\ref{Sch2}) and letting $p:=2r+p',\,(r\in\ZZ^g/n\ZZ^g,p'\in\ZZ^g/2\ZZ^g)$,
one has by virtue of (\ref{d3a}):
\be
\theta_{[a_1;b_1]}(u_1|\tau)\theta_{[a_2;b_2]}(u_2|\tau)=
n^{-g}\sum_{p'\in\ZZ^g/2\ZZ^g}\theta_{[\frac{a_1+a_2+p'}{2};b_1+b_2]}(u_1+u_2|2\tau)
\\ \times
\sum_{q\in\ZZ^g/n\ZZ^g}e^{-\frac{\pi i}{n}\<a_1-a_2+p',q\>}
\theta_{[\frac{a_1-a_2+p'}{2};b_1-b_2+\frac{q}{n}]}(u_1+u_2|2\tau)
\sum_{r\in\ZZ^g/n\ZZ^g}e^{-\frac{2\pi i}{n}\<r,q\>}.
\ee
Using the identity (\ref{sr3}), one arrives  to the standard binary identities \cite{M1}, \cite{D}:
\bel{bi1}
\theta_{[a_1;b_1]}(u_1|\tau)\theta_{[a_2;b_2]}(u_2|\tau)=\\
\sum_{p\in\ZZ^g/2\ZZ^g}
\theta_{[\frac{a_1+a_2+p}{2};b_1+b_2]}(u_1+u_2|2\tau)
\theta_{[\frac{a_1-a_2+p}{2};b_1-b_2]}(u_1-u_2|2\tau).
\ee
In the rest of the paper, we shall deal only with the standard binary identities (\ref{bi1}).
\begin{cor}
The inverse binary relations are:
\bel{bi2}
\theta_{[a_1;b_1]}(u_1+u_2|2\tau)
\theta_{[a_2;b_2]}(u_1-u_2|2\tau)\\=
2^{-g}\sum_{p\in\ZZ^g/2\ZZ^g}e^{-2\pi i\<a_1,p\>}
\theta_{[a_1+a_2;\frac{b_1+b_2+p}{2}]}(u_1|\tau)
\theta_{[a_1-a_2;\frac{b_1-b_2+p}{2}]}(u_2|\tau).
\ee
\end{cor}
{\bf Proof.} In (\ref{bi1}), perform the
shifts $b_k\to b_k+\frac{q}{2}$, where by definition $q\in\ZZ^g$. Due to (\ref{d3a}), the product of theta functions in the right hand side of (\ref{bi1}) acquires the form $e^{\pi i\<a_1+a_2+p,q\>}$ $\theta_{[\frac{a_1+a_2+p}{2};b_1+b_2]}(u_1+u_2|2\tau)
\theta_{[\frac{a_1-a_2+p}{2};b_1-b_2]}(u_1-u_2|2\tau)$. Thus we have:
\bel{bi3}
e^{-\pi i\<a_1+a_2,q\>}
\theta_{[a_1;b_1+\frac{q}{2}]}(u_1|\tau)\theta_{[a_2;b_2+\frac{q}{2}]}(u_2|\tau)\\=
\sum_{p\in\ZZ^g/2\ZZ^g}e^{\pi i\<p,q\>}
\theta_{[\frac{a_1+a_2+p}{2};b_1+b_2]}(u_1+u_2|2\tau)
\theta_{[\frac{a_1-a_2+p}{2};b_1-b_2]}(u_1-u_2|2\tau).
\ee
Note that the right hand side of (\ref{bi3}) depends on the vector $q\in\ZZ^g$
in a very simple manner.
In particular, performing the summation over $q\in\ZZ^g/2\ZZ^g$ in (\ref{bi3}) and using the  identity
\bel{bi4}
\sum_{q\in\ZZ^g/2\ZZ^g}e^{\pi i\<p,q\>}=2^g\delta_{p,0},\ \ \ \ p\in\ZZ^g/2\ZZ^g,
\ee
one arrives to (\ref{bi2}). \ws
\begin{rem}
The derivation of Jacobi and Riemann identities below is essentially the same as derivation of (\ref{bi2}): taking the appropriate product of theta functions, the only thing is to perform the relevant shifts of characteristics $[a;b]$ and take into account the simple properties (\ref{d3}). Then all identities arise by virtue of identity (\ref{bi4}).
\end{rem}

\section{Jacobi identities}
In one-dimensional case ($g=1$), the four-term identities of order 4 were essentially
obtained by Jacobi \cite[p. 507]{J2} (see also \cite[pp. 468, 488]{WW}). For example, there is the identity
\bel{J1}
\prod_{k=1}^4\theta_{[\frac{1}{2};\frac{1}{2}]}(u_k|\tau)+
\prod_{k=1}^4\theta_{[\frac{1}{2};0]}(u_k|\tau)=
\prod_{k=1}^4\theta_{[\frac{1}{2};\frac{1}{2}]}(u_k'|\tau)+
\prod_{k=1}^4\theta_{[\frac{1}{2};0]}(u_k'|\tau),
\ee
where the primed variables $u_k',\,k=1,2,3,4$ are defined by (\ref{in3}).
The problem is to generalize the relations of type (\ref{J1}) to multidimensional case.
In this section we derive the identities (\ref{in5}) and (\ref{in6}) from (\ref{bi1}) and (\ref{bi2}), respectively.

The derivation of (\ref{in5}) is as follows. Shifting in  (\ref{bi1}) $b_k\to b_k+\frac{q}{2},\,q\in\ZZ^g/2\ZZ^g$ and taking into account (\ref{d3a}), one has:
\bel{J2}
\prod_{k=1}^4e^{-\pi i\<a_k,q\>}\theta_{[a_k;b_k+\frac{q}{2}]}(u_k|\tau)\\=
\sum_{p,p'\in\ZZ^g/2\ZZ^g}e^{\pi i\<p+p',q\>}
\theta_{[\frac{a_1+a_2+p}{2};b_1+b_2]}(u_1+u_2|2\tau)
\theta_{[\frac{a_1-a_2+p}{2};b_1-b_2]}(u_1-u_2|2\tau)\\ \times\,
\theta_{[\frac{a_3+a_4+p'}{2};b_3+b_4]}(u_3+u_4|2\tau)
\theta_{[\frac{a_3-a_4+p'}{2};b_3-b_4]}(u_3-u_4|2\tau).
\ee
Note that the shift $b_k\to b_k+\frac{q}{2}$ implies the transformation
$b_k'\to b_k'+\frac{q}{2}$, where the Whittaker-Watson dual vectors $b_k'$ are defined in accordance with (\ref{in3}). Therefore, by virtue of (\ref{d3}) one has:
\bel{J3}
\prod_{k=1}^4e^{-\pi i\<a_k',q\>}\theta_{[a_k';b_k'+\frac{q}{2}]}(u_k'|\tau)\\=
\sum_{p,p'\in\ZZ^g/2\ZZ^g}e^{\pi i\<p+p',q\>}
\theta_{[\frac{a_1+a_2+p'}{2};b_1+b_2]}(u_1+u_2|2\tau)
\theta_{[\frac{a_1-a_2-p}{2};b_1-b_2]}(u_1-u_2|2\tau)\\ \times\,
\theta_{[\frac{a_3+a_4+p}{2};b_3+b_4]}(u_3+u_4|2\tau)
\theta_{[\frac{a_3-a_4-p'}{2};b_3-b_4]}(u_3-u_4|2\tau).
\ee
Performing the summation over $q\in\ZZ^g/2\ZZ^g$
in (\ref{J2}), (\ref{J3}) and using the identity (\ref{bi4}), one arrives to
the following expressions:
\bel{J4}
\sum_{q\in\ZZ^g/2\ZZ^g}\prod_{k=1}^4
e^{-\pi i\<a_k,q\>}\theta_{[a_k;b_k+\frac{q}{2}]}(u_k|\tau)\\ =2^g
\sum_{p\in\ZZ^g/2\ZZ^g}
\theta_{[\frac{a_1+a_2+p}{2};b_1+b_2]}(u_1+u_2|2\tau)
\theta_{[\frac{a_1-a_2+p}{2};b_1-b_2]}(u_1-u_2|2\tau)\\ \times\,
\theta_{[\frac{a_3+a_4-p}{2};b_3+b_4]}(u_3+u_4|2\tau)
\theta_{[\frac{a_3-a_4-p}{2};b_3-b_4]}(u_3-u_4|2\tau),
\ee
\bel{J5}
\sum_{q\in\ZZ^g/2\ZZ^g}\prod_{k=1}^4
e^{-\pi i\<a_k',q\>}\theta_{[a_k';b_k'+\frac{q}{2}]}(u_k'|\tau)\\ =2^g
\sum_{p\in\ZZ^g/2\ZZ^g}
\theta_{[\frac{a_1+a_2-p}{2};b_1+b_2]}(u_1+u_2|2\tau)
\theta_{[\frac{a_1-a_2-p}{2};b_1-b_2]}(u_1-u_2|2\tau)\\ \times\,
\theta_{[\frac{a_3+a_4+p}{2};b_3+b_4]}(u_3+u_4|2\tau)
\theta_{[\frac{a_3-a_4+p}{2};b_3-b_4]}(u_3-u_4|2\tau).
\ee
Finally, the right hand sides of (\ref{J4}) and (\ref{J5}) coincide by virtue of (\ref{d3a}) and identities (\ref{in5}) do hold.

Quite similarly, one can derive identities (\ref{in6}) staring with the Schr\"oter relations (\ref{bi2}) written in the form
\bel{J6}
\theta_{[a_1;b_1]}(u_1|\tau)\theta_{[a_2;b_2]}(u_2|\tau)\\=
2^{-g}\sum_{p\in\ZZ^g/2\ZZ^g}e^{-2\pi i\<a_1,p\>}
\theta_{[a_1+a_2;\frac{b_1+b_2+p}{2}]}(\frac{u_1+u_2}{2}|\frac{\tau}{2})
\theta_{[a_1-a_2;\frac{b_1-b_2+p}{2}]}(\frac{u_1-u_2}{2}|\frac{\tau}{2}).
\ee
Hence, the relations (\ref{in5}), (\ref{in6}) are simple algebraic consequences of
the binary
Schr\"oter identities. Below we prove the equivalence of (\ref{in5}) and (\ref{in6}). \ws

\bigskip
For completeness, we represent the Jacobi identities in a slightly different form.
In accordance with original work \cite[p. 503]{J2}, introduce the Jacobi dual variables as follows:
\bel{J7}
\tilde x_1={\textstyle\frac{1}{2}}(x_1+x_2+x_3+x_4),\\
\tilde x_2={\textstyle\frac{1}{2}}(x_1+x_2-x_3-x_4),\\
\tilde x_3={\textstyle\frac{1}{2}}(x_1-x_2+x_3-x_4),\\
\tilde x_4={\textstyle\frac{1}{2}}(x_1-x_2-x_3+x_4).
\ee
These variables have been used, for example, by  Prym \cite[p. 6]{Pr} in his monograph on multidimensional Riemann identities. The same notations are used also in \cite[p. 212]{M1}, \cite[p. 102]{M2}.

Evidently, the change $x_1\to -x_1$ in (\ref{in3}) leads to the transformation
of the Whittaker-Watson dual variables to the Jacobi ones:
$x_1'\to \tilde x_1,x_2'\to-\tilde x_2,x_3'\to-\tilde x_3,x_4'\to-\tilde x_4$. As a consequence, it is easy to see that the Jacobi identities can be written the form:
\bs{J8}
\ba
\sum_{p\in\ZZ^g/2\ZZ^g}
\prod_{k=1}^4\theta_{[a_k+\frac{p}{2};b_k]}(u_k|\tau)&=
\sum_{p\in\ZZ^g/2\ZZ^g}
\prod_{k=1}^4\theta_{[\tilde a_k+\frac{p}{2};\tilde b_k]}(\tilde u_k|\tau),\label{J8a}\\
\sum_{q\in\ZZ^g/2\ZZ^g}\prod_{k=1}^4
e^{-\pi i\<a_k,q\>}\theta_{[a_k;b_k+\frac{q}{2}]}(u_k|\tau)&=
\sum_{q\in\ZZ^g/2\ZZ^g}\prod_{k=1}^4
e^{-\pi i\<\tilde a_k,q\>}\theta_{[\tilde a_k;\tilde b_k+\frac{q}{2}]}(\tilde u_k|\tau)
\label{J8b}.
\end{align}
\es
Indeed, one can change $(a_1,b_1,u_1)\to(-a_1,-b_1,-u_1)$ in (\ref{in5}). Then, by virtue of (\ref{d3b}) and (\ref{d3a}), the identity (\ref{in5}) gives (\ref{J8a}). Similarly, under the same change, the identity (\ref{J8b}) follows from (\ref{in6}).

\section{Riemann identities}
One can derive the Riemann identities (\ref{in7}) from the Jacobi's ones (\ref{in5}) or (\ref{in6}). For example, perform the following shifts in (\ref{in6}):
$b_1\to b_1-\frac{q}{2},\,
b_2\to b_2+\frac{q}{2},\,b_3\to b_3+\frac{q}{2},\,b_4\to b_4+\frac{q}{2}$, where $q\in\ZZ^g/2\ZZ^g$. According to
the definition of the Whittaker-Watson dual variables (\ref{in3}), one has the corresponding transformations $b_1'\to b_1'+q,\,
b_2'\to b_2',\,b_3'\to b_3',\,b_4'\to b_4'$. Taking into account (\ref{d3a}), the Jacobi identity (\ref{in6}) acquires the form:
\bel{R1}
e^{-2\pi i\<a_1+a_1',q\>}\sum_{p\in\ZZ^g/2\ZZ^g}e^{-\pi i\<p,q\>}
\prod_{k=1}^4\theta_{[a_k+\frac{p}{2};b_k+\frac{q}{2}]}(u_k|\tau)\\=
\sum_{p\in\ZZ^g/2\ZZ^g}e^{\pi i\<p,q\>}
\prod_{k=1}^4\theta_{[a_k'+\frac{p}{2};b_k']}(u_k'|\tau).
\ee
Performing the summation over $q\in\ZZ^g/2\ZZ^g$
and using identity (\ref{bi4}), one arrives at the relation
\bel{R2}
2^g\prod_{k=1}^4\theta_{[a_k';b_k']}(u_k'|\tau)=
\sum_{p,q\in\ZZ^g/2\ZZ^g}e^{-\pi i\<p+2a_1+2a_1',q\>}
\prod_{k=1}^4\theta_{[a_k+\frac{p}{2};b_k+\frac{q}{2}]}(u_k|\tau).
\ee
Noting that in accordance with definition (\ref{in3}) $2a_1+2a_1'=a_1+a_2+a_3+a_4$, the relation (\ref{R2}) gives the Riemann identity (\ref{in7}). Quite similarly, the same identity can be derived from (\ref{in5}) performing the shifts
$a_1\to a_1-\frac{p}{2},\,
a_2\to a_2+\frac{p}{2},\,a_3\to a_3+\frac{p}{2},\,a_4\to a_4+\frac{p}{2}$, where $p\in\ZZ^g/2\ZZ^g$.

It is easy to see that the connection between $u_k$ and $u_k'$ (\ref{in3}) is a reciprocal one. Hence, the inverse Riemann identities are:
\bel{R3}
\hspace{-0.3cm}
\prod_{k=1}^4\theta_{[a_k;b_k]}(u_k|\tau)=2^{-g}
\sum_{p,q\in\ZZ^g/2\ZZ^g}e^{-\pi i\<p,q\>}\prod_{k=1}^4
e^{-\pi i\<a_k',q\>}\theta_{[a_k'+\frac{p}{2};b_k'+\frac{q}{2}]}(u_k'|\tau).
\ee

We should remark that in terms of Jacobi variables (\ref{J7}), the Riemann identities read
(see, for example, \cite[p. 166]{T}):
\bel{R4}
\prod_{k=1}^4\theta_{[\tilde a_k;\tilde b_k]}(\tilde u_k|\tau)=
2^{-g}\sum_{p,q\in\ZZ^g/2\ZZ^g}\prod_{k=1}^4e^{-\pi i\<a_k,q\>}
\theta_{[a_k+\frac{p}{2};b_k+\frac{q}{2}]}(u_k|\tau).
\ee
Indeed, (\ref{R4}) can be easily obtained from (\ref{in7}) by the change
$(a_1,b_1,u_1)\to(-a_1,-b_1,-u_1)$.
Finally, it is easy to see that the connection between $u_k$ and $\tilde u_k$ (\ref{J7}) is a reciprocal one. Hence, the inverse formula to (\ref{R4}) is:
\bel{R5}
\prod_{k=1}^4\theta_{[a_k;b_k]}(u_k|\tau)=
2^{-g}\sum_{p,q\in\ZZ^g/2\ZZ^g}\prod_{k=1}^4e^{-\pi i\<\tilde a_k,q\>}
\theta_{[\tilde a_k+\frac{p}{2};\tilde b_k+\frac{q}{2}]}(\tilde u_k|\tau).
\ee

\section{Equivalence of the Jacobi and Riemann identities}
In the previous section, we have obtained the
Riemann identities (\ref{in7}) from the
Jacobi identities of the 'second kind' (\ref{in6}). It is easy to show
that (\ref{in7}) implies both (\ref{in5}) and (\ref{in6}).
For example, perform the shifts $a_k\to a_k+r\delta_{k,1},\,k=1,2,3,4$,
where $r\in\ZZ^g/2\ZZ^g$. Then $a_k'\to a_k'+(-1)^{\delta_{k,1}}\frac{r}{2}$, and
by virtue of (\ref{d3a}) the identity (\ref{in7}) acquires the form:
\bel{eq1}
\hspace{-0.5cm}
\prod_{k=1}^4\theta_{[a_k'+\frac{r}{2};b_k']}(u_k'|\tau)=2^{-g}
\sum_{p,q\in\ZZ^g/2\ZZ^g}e^{-\pi i\<p,q\>-\pi i\<r,q\>}\prod_{k=1}^4
e^{-\pi i\<a_k,q\>}\theta_{[a_k+\frac{p}{2};b_k+\frac{q}{2}]}(u_k|\tau).
\ee
Performing in (\ref{eq1}) the summation over $r\in\ZZ^g/2\ZZ^g$ and using (\ref{bi4}), one arrives to the Jacobi identity (\ref{in6}). Similarly, performing the shifts
$b_k\to b_k+r\delta_{k,1},\,k=1,2,3,4$ in (\ref{in7}), one arrives to (\ref{in5}). Thus the Jacobi and Riemann identities are equivalent.
\ws

\bigskip
Finally, one can represent the Riemann identities in slightly different form which relates the products of theta functions in Whittaker-Watson and Jacobi variables. Changing in the left hand side of (\ref{R5}) $(a_k,b_k,u_k)\to(a_k',b_k',u_k')$, one obtains the following transformations of all variables in the right hand side of (\ref{R5}):
$\tilde x_1\to\tilde x_1,\tilde x_2\to-\tilde x_2, \tilde x_3\to-\tilde x_3,
\tilde x_4\to-\tilde x_4$.
Therefore, by virtue of (\ref{d3}), one can rewrite (\ref{R5}) as follows:
\bel{R6}
\prod_{k=1}^4\theta_{[a_k';b_k']}(u_k'|\tau)=
2^{-g}\sum_{p,q\in\ZZ^g/2\ZZ^g}e^{-\pi i\<p,q\>}
\prod_{k=1}^4e^{-\pi i\<\tilde a_k,q\>}
\theta_{[\tilde a_k+\frac{p}{2};\tilde b_k+\frac{q}{2}]}(\tilde u_k|\tau).
\ee
As a trivial corollary of (\ref{R5}) and (\ref{R6}), the following identities hold:
\bel{R7}
\prod_{k=1}^4\theta_{[a_k;b_k]}(u_k|\tau)-\prod_{k=1}^4\theta_{[a_k';b_k']}(u_k'|\tau)\\
=2^{-g}\sum_{p,q\in\ZZ^g/2\ZZ^g}\big(1-e^{-\pi i\<p,q\>}\big)
\prod_{k=1}^4e^{-\pi i\<\tilde a_k,q\>}
\theta_{[\tilde a_k+\frac{p}{2};\tilde b_k+\frac{q}{2}]}(\tilde u_k|\tau).
\ee
Certainly, the relations (\ref{R7}) are equivalent to the Jacobi identities (\ref{in5}), (\ref{in6}).
Indeed, perform the shifts
$a_k\to a_k+\frac{r}{2},\,r\in\ZZ^g/2\ZZ^g$.
Then in accordance with (\ref{in3}), (\ref{J7}), $a_k'\to a_k'+\frac{r}{2}$,
$\tilde a_k\to\tilde a_k+r\delta_{k,1}$, and by virtue of (\ref{d3a}) the relation (\ref{J7}) acquires the form
\bel{R8}
\prod_{k=1}^4\theta_{[a_k+\frac{r}{2};b_k]}(u_k|\tau)-
\prod_{k=1}^4\theta_{[a_k'+\frac{r}{2};b_k']}(u_k'|\tau)\\
=2^{-g}\sum_{p,q\in\ZZ^g/2\ZZ^g}\big(1-e^{-\pi i\<p,q\>}\big)
e^{-\pi i\<r,q\>}
\prod_{k=1}^4e^{-\pi i\<\tilde a_k,q\>}
\theta_{[\tilde a_k+\frac{p}{2};\tilde b_k+\frac{q}{2}]}(\tilde u_k|\tau).
\ee
Performing summation over $r\in\ZZ^g/2\ZZ^g$ it is easy to see that the right hand side of (\ref{R8}) vanishes by virtue of (\ref{bi4}). Hence, one arrives to the Jacobi identities (\ref{in6}).
\begin{ex}{\rm
In one-dimensional case ($g=1$) the general formulas (\ref{R7}) read:
\bel{R9}
\prod_{k=1}^4\theta_{[a_k;b_k]}(u_k|\tau)-\prod_{k=1}^4\theta_{[a_k';b_k']}(u_k'|\tau)=
e^{-2\pi i a_1}\prod_{k=1}^4
\theta_{[\tilde a_k+\frac{1}{2};\tilde b_k+\frac{1}{2}]}(\tilde u_k|\tau),
\ee
which describe all Weierstrass addition formulas listed in {\rm\cite{KZ}}. In particular, letting $a_k=b_k=\frac{1}{2}$ and using the standard notation
$\theta_1(u|\tau):=-\theta_{[\frac{1}{2};\frac{1}{2}]}(u|\tau)$, one has:
\bel{ex1}
\prod_{k=1}^4\theta_1(u_k)-\prod_{k=1}^4\theta_1(u_k')=
\prod_{k=1}^4\theta_1(\tilde u_k).
\ee
To rewrite (\ref{ex1}) in a
more explicit form, set $u_1=w_1+w_2,u_2=w_1-w_2,u_3=w_3+w_4,u_4=w_3-w_4$.
Then the relation (\ref{ex1}) has the form
\bel{ex2}
\theta_1(w_1+w_2)\theta_1(w_1-w_2)\theta_1(w_3+w_4)\theta_1(w_3-w_4)\\-
\theta_1(w_3-w_2)\theta_1(w_3+w_2)\theta_1(w_1-w_4)\theta_1(w_1+w_4)\\=
\theta_1(w_1+w_3)\theta_1(w_1-w_3)\theta_1(w_2+w_4)\theta_1(w_2-w_4).
\ee
This is the famous Weierstrass identity for the function $\sigma(u)$ \cite{We1}.
In preface of his paper, Weierstrass pointed out that identity (\ref{ex2}) is essentially different from the Jacobi ones since the latter contain two or more theta functions.}
\footnote{
\;"Diese Gleichung ist wesentlich anderer Art als die von Jacobi entdeckten, auf S. 507 des ersten Bandes der "Gesammelten Werke" vollstandig augestellten Relationen unter Producten von je vier $\theta$-Functionen; sie enthalt nur eine Function, wahrend in jeder der Jacobi'sehen Gleichungen, die sich ubrigens aus ihr ableiten lassen, zwei oder mehrere $\theta$-Functionen vorkommen". \cite[p. 155]{We2}
}
\end{ex}

Nevertheless it turns out \cite{K}, \cite{KZ} that in the
one-dimensional case the Jacobi and Weierstrass identities are equivalent.
In particular, the addition formula (\ref{ex2}) is a simple corollary of Jacobi identity (\ref{J1}) (see \cite{KZ} for detail).

\section{Weierstrass identities}
The problem raised by Weierstrass \cite{We2} is to find the identity for the single
multidimensional theta function 
thus obtaining the direct generalization of identity (\ref{ex2}). From this point of view, the naive multidimensional identities (\ref{R7}) do not fit the scheme.

Let $[a;b]$ be half-periods, i.e the components of $g$-dimensional vectors $a$ and $b$ are 0 or $\frac{1}{2}$.
\begin{prop}\label{pw}
{\rm\cite{We2}} Let $\theta(u)$ be any $g$-dimensional odd theta function and
$w_i,\,i=1,\ldots,2^g+2$ be arbitrary variables. Then the Pfaffian of $(2^g+2)\times(2^g+2)$ skew-symmetric matrix $A=||\theta(w_i+w_j)\theta(w_i-w_j)||_{i,j=1}^{2^g+2}$
identically vanishes:
\bel{mp1}
{\rm Pf}||\theta(w_i+w_j)\theta(w_i-w_j)||_{i,j=1}^{2^g+2}=0.
\ee
\end{prop}
The proof of Proposition \ref{pw} is a simple corollary of the following Lemma:
\begin{lem}\label{lem}
Let $\theta(u)$ be an arbitrary odd $g$-dimensional theta function. Then
\bel{f1}
\theta(w_1+w_2)\theta(w_1-w_2)=\sum_{k=1}^{2^{g-1}}
\Big\{A_k(w_1)B_k(w_2)-B_k(w_1)A_k(w_2)\Big\},
\ee
where $A_k(w),B_k(w),k=1,\ldots, 2^{g-1}$ are appropriate theta functions of type
$\theta_{[c;0]}(2w|2\tau)$, $c\in\frac{1}{2}\ZZ^g/2\ZZ^g$.
\end{lem}
{\bf Proof} of Lemma \ref{lem}.
Consider the multidimensional Schr\"oter identities (\ref{bi1}). For any half-periods $[a;b]$ one has the particular relations
\bel{ps}
\theta_{[a;b]}(w_1+w_2|\tau)\theta_{[a;b]}(w_1-w_2|\tau)\\=
e^{4\pi i\<a,b\>}\sum_{p\in\ZZ^g/2\ZZ^g}e^{2\pi i\<p,b\>}
\theta_{[a+\frac{p}{2};0]}(2w_1|2\tau)
\theta_{[\frac{p}{2};0]}(2w_2|2\tau)
\ee
which hold by virtue of (\ref{d3a}). The sum in (\ref{ps}) contains $2^g$ terms, and for any given expression
$e^{2\pi i\<p,b\>}\theta_{[a+\frac{p}{2};0]}(2w_1|2\tau)\theta_{[\frac{p}{2};0]}(2w_2|2\tau)$ there is the corresponding term
$$e^{2\pi i\<p+2a,b\>}\theta_{[a+\frac{p+2a}{2};0]}(2w_1|2\tau)
\theta_{[\frac{p+2a}{2};0]}(2w_2|2\tau)$$
$$=e^{4\pi i\<a,b\>}\cdot e^{2\pi i\<p,b\>}
\theta_{[\frac{p}{2};0]}(2w_1|2\tau)\theta_{[a+\frac{p}{2};0]}(2w_2|2\tau)
$$
in the sum (\ref{ps}). Hence for the odd periods $4\<a,b\>=1\,({\rm mod}\;2)$ one has the structure (\ref{f1}) and Lemma \ref{lem} is proved. \ws

Now the proof of Proposition \ref{pw} is very simple. It is easy to see that the determinant of $(2n)\times(2n)$ skew-symmetric matrix $||a_{i,j}||$ identically vanishes provided the structure
$a_{i,j}=\sum_{k=1}^{n-1}\{A_k(w_1)B_k(w_2)-B_k(w_1)A_k(w_2)\}$. \ws

\begin{ex}
The Pfaffian of $4\times4$ skew-symmetric matrix $||a_{i,j}||$ is:
\be
{\rm Pf}\,||a_{i,j}||=a_{1,2}a_{3,4}-a_{1,3}a_{2,4}+a_{1,4}a_{2,3}.
\ee
Letting $a_{i,j}=\theta_1(w_i+w_j)\theta_1(w_i-w_j)$, one arrives to the Weierstrass identity (\ref{ex2}) in the case $g=1$.
\end{ex}


\bibliographystyle{12}
\bibliographystyle{amsalpha}

\end{document}